\documentclass{elsart}
\usepackage{amssymb}
\usepackage{mathrsfs}
\usepackage{graphicx}
\usepackage{stmaryrd}
\usepackage{dsfont}
\usepackage{algorithm,algorithmic}

\newtheorem{rque}{Remark}[section]

\newcommand{\R}{\mathbb{R}}
\newcommand{\N}{\mathbb{N}}
\newcommand{\Z}{\mathbb{Z}}
\newcommand{\supp}{\mbox{\rm{supp }}}
\journal{}
\begin{document}
\begin{frontmatter}
\title{A Construction of Biorthogonal Wavelets With a Compact Operator}
\author[LAMA]{Mehmet Ersoy}
\ead{Mehmet.Ersoy@univ-savoie.fr}
\address[LAMA]{Universit\'{e} de Savoie,
Laboratoire de Math\'{e}matiques, 73376 Le Bourget-du-Lac Cedex, France}
\begin{abstract}
We present a construction of  biorthogonal wavelets using a compact operator  which allows to preserve or increase some properties:
regularity/vanishing moments, parity, compact supported. We build then a simple algorithm which computes new filters.
\end{abstract}
\end{frontmatter}
\section{Introduction}
The aim of this paper is to present a construction of   biorthogonal wavelets with some properties of vanishing moments/regularity
which could be a powerful tool in wavelet partial differential equations theory ;
in particular in the case of the wavelet-Galerkin method (WG)(see e.g \cite{AWQW94})
or the so-called hybrid method (\cite{FS97,PL96}). The property of regularity/vanishing moment is more important in wavelets theory.
To this end, we use an iterative process to get new wavelets.
The new wavelets have $s$ vanishing moments while the biorthogonal ones are more regular.
Moreover, this construction could be a way to recover all boxed spline wavelets.
As in most works, we privilege compactly supported wavelets which are very attractive for such a construction.\\
We present in the following section a way to construct these wavelets using a simple compact operator. To this end, we choose initially
biorthogonal wavelets (orthonormal wavelets can also be choose) and we define them from the so-called transfer functions. We show how
to build easily theirs discrete filters in practice.
In what follows, we denote by: BMRA biorthogonal multi-resolution analysis, OMRA orthogonal multi-resolution analysis, ($\phi,\,\psi$)
initial scaling function-wavelet pair with ($\widetilde\phi,\,\widetilde\psi$) the dual pair, ($\Phi^{(s)},\,\Psi^{(s)}$) the constructed
scaling function-wavelet pair at order $s$ where $s$ design the vanishing moments (or degree of approximation)  with ($\widetilde\Phi^{(s)},\,\widetilde\Psi^{(s)}$)
the dual pair where $s$ denotes the order of regularity  and $\langle \cdot\,,\, \cdot\rangle$
denotes the usual scalar product of $L^2(\R)$. We denote also by (TSR) the classical two scale relation
\section{Main results}
Our motivation is to construct a high order method to elevate
the degree of vanishing moment and regularity of any given biorthogonal wavelet
in order to get powerful wavelets for wavelet partial differential equations theory.
For this purpose, let us introduce the compact operator $T$ of kernel $G(x,y)=\mathds{1}_{[x-1,x]}(y)$
where its adjoint is defined by $G^{*}(x,y)= \mathds{1}_{[x,x+1]}(y)$.
The operator $T$ acting $s$ times on the initial wavelet is an elegant way  to build a new pair of biorthogonal wavelets.
The process is called  \emph{an elevation at order $s$}. To start the construction, we choose initially biorthogonal wavelets.
\begin{thm}\label{thm_principal}
Let $\phi$,\,$\tilde{\phi}$ be two scaling functions  of a BMRA. Let $T$ be the operator described above and
$T^{*}$ its adjoint. Let $s\in\N$ be the order of elevation and  let us define $\Phi^{(s)}$ and $\tilde{\Phi}^{(s)}$ as
$$ \Phi^{(s)}:= (T)^{s} \phi,\quad  \tilde{\Phi}^{(s)}:= \left(T^{*}\right)^{-s}\tilde{\phi} $$ where $(X)^s$
designs the action of $X $ $s$
times, and by definition $X^{0}:=Id$.\\
Then, the scaling functions  $\Phi^{(s)}$ and $\tilde{\Phi}^{(s)}$ generates a BMRA. Moreover,
$$P_0^{(s)}(w) = \frac{1}{2^s} \frac{S_s(2w)}{S_s(w)} m_0(w),\,
 \tilde{P}_0^{(s)}(w) = 2^s \frac{\overline{S_s(w)}}{\overline{S_s(2w)}} \tilde{m}_0(w) $$
$$P_1^{(s)}(w) = 2^s \frac{m_1(w)}{S_s(w)},\,\tilde{P}_1^{(s)}(w) = \frac{1}{2^s} \overline{S_s(w)}
\tilde{m}_1(w)$$  where
$m_0$, $\tilde{m}_0$, $P_0^{(s)}$, $\tilde{P}_0^{(s)}$, $m_1,\,\tilde{m}_1$, $P_1^{(s)}$, $\tilde{P}_1^{(s)}$
  denote the transfer functions of
$\phi$, $\tilde{\phi}$, $\Phi^{(s)}$, $\tilde{\Phi}^{(s)}$, $\psi$, $\tilde{\psi}$, $\Psi^{(s)}$, $\tilde{\Psi}^{(s)}$
and  $S_s(w) =(1-e^{-iw})^s.$
\end{thm}
To prove Theorem \ref{thm_principal}, let us recall the following results:
\begin{lem}\label{gamma} Let $\phi$ be a scaling function then
\begin{description}
 \item[$(i)$]  $\Gamma >0$ where $\Gamma(w)=\displaystyle \sum_{k\in\mathbb{Z}} |\hat{\Phi}(w+2k\pi)|^2$.
 \item[$(ii)$] The set $\{\Phi(.-k),\,k\in\Z\}$ forms a Riesz basis.
\end{description}
\end{lem}
\begin{lem}\label{Yen_lemma} Let $(\phi,\psi)$ be the scaling function and wavelet of an OMRA. We suppose
$\phi$, $\psi$ differentiable. Then the following equality holds
$$\langle \psi'(x),\psi'(x-k)\rangle + \langle \phi'(x),\phi'(x-k)\rangle = 4 \langle \phi'(x),\phi'(x-2k)\rangle\,.$$
\end{lem}
In order to complete the proof let us do the following remark.
\begin{rque}\rm
We have:
\begin{equation}\label{regular} \Psi(x)=4 \int_{-\infty}^x \psi(t)dt,\,\,\,
4 \int_{-\infty}^x \tilde\Psi(t)dt=\psi(x),\,\,\, \forall x\in\mathbb{R}.\end{equation}
\end{rque}
\noindent \textbf{Proof of theorem \ref{thm_principal}: } For $s=0$, the results holds since we get the classical definition of the
transfer functions. For $s=1$, see e.g. \cite{L92}. For $s>1$, let us denote to simplify notations $\phi=\Phi^{(s)}$ and $\Phi=\Phi^{(s+1)}$.\\
We show firstly that the set $\{\Phi(\cdot-k),\,k\in\Z\}$ cannot be
orthonormal and thus there exists a biorthogonal function denoted by $\tilde{\Phi}$.
To this end, we proceed as follows:
\begin{description}
 \item[Step $1$. ] $\Phi$ satisfies a the (TSR). Indeed, since
$$\Phi(x)= T\phi(x) = \int_{-\infty}^{x}
\left(\phi(t)-\phi(t-1)\right)dt $$ and using the (TSR) on $\phi$.
Furthermore, when $\supp\phi^{(s)}= [0,2p-1+s]$, we obtain the
discrete filters of $\Phi$ defined as:
$$H_k=\displaystyle\frac{h_{k-1}+h_k}{2},\,k=0\ldots 2p \textrm{
with } h_{-1}=h_{2p}=0.$$
\item[Step $2$. ] Now, it remains to show  that
$(\Phi(\cdot-k))_{k\in\Z}$ forms a Riesz basis. For this purpose, we use Lemma \ref{gamma} to find two
constants $0<A \leq B$ such that $$\forall
w\in\mathbb{R},\,\displaystyle A\leq \sum_{k\in\mathbb{Z}}
|\hat{\Phi}(w+2k\pi)|^2 \leq B.$$
The fact that $\Gamma$ is $2\pi$-periodic continuous
with compacity property give us the constant $B$ whereas $A$ is obtained since $\Gamma>0$.
\item[Step $3$. ] To conclude, we  prove that $(\Phi(\cdot-k))_{k\in\Z}$
never forms an orthonormal family but it satisfy $
\displaystyle \langle \Phi(\cdot-k),\Phi(\cdot)\rangle=0 $ for almost all integer $ k$
which ensures that the family $\Phi$ cannot be orthonormal. This is done by setting
$\tilde{\phi}=T^{\ast} \tilde\Phi$  and using Property (\ref{regular}), Lemma \ref{Yen_lemma}.
The transfer functions are obtained by a straightforward computation.\\
\end{description}
When $\phi,\,\tilde\phi$ are compactly supported, we have:
\begin{thm}\label{thm_principal2} Under the hypothesis of theorem \ref{thm_principal}, if we suppose that $\phi$ and $\tilde{\phi}$ are compactly supported where $\supp\phi=[0,2p-1]$ and $\supp\tilde{\phi}=[0,2\tilde p-1]$ then
writing $H_k^{(s)}$, $\tilde{H}_k^{(s)}$ the discrete filters of  $\Phi^{(s)}$, $\tilde{\Phi}^{(s)}$  (resp.)
and  $h_k^{(s)}$, $\tilde{h}_k^{(s)}$ the discrete filters of $\phi^{(s)}$, $\tilde{\phi}^{(s)}$ we have,
 $$\begin{array}{lll}
 \supp\Phi^{(s)} &=& [0,2p-1+s] \\
 \supp\tilde{\Phi}^{(s)} &=& [0,2\tilde{p}-1-s]\end{array}$$
with the following formulas to compute the new filters from the older
\begin{equation}
\begin{array}{lll}\label{formul1}
 H_k^{(s)} = \displaystyle \frac{1}{2^s} \sum_{l=0}^{s}
\left( \begin{array}{c} s\\l\end{array}\right)h_{k-l}, &  k=0\ldots 2p-1+s
\end{array}
\end{equation}
\begin{equation}\label{formul2}
\begin{array}{lll}
\displaystyle \sum_{l=0}^{s} \left( \begin{array}{c} s\\l\end{array}\right) \tilde{H}^{(s)}_{k-l} =2^s
\tilde{h}_k,& k=0\ldots 2\tilde p-1-s,\,\tilde p\geq
\displaystyle\frac{s+1}{2}
\end{array}
\end{equation}
\end{thm}
\noindent\textbf{Proof of Theorem \ref{thm_principal2}: }
We have already proved the result for $s=1$ (\textbf{Step $1$.} of the proof of Theorem \ref{thm_principal}).
Writing $H_k^{0}:=h_k$, we have
\begin{description}
\item[ for $s = 1$: ] $\displaystyle H_k^{(1)} =
\frac{H^{(0)}_{k-1} + H^{(0)}_{k}}{2}$ \- and $\displaystyle
\Phi^{(1)}$ admits an approximation of degree $p+1$ and
$\supp\displaystyle \Phi^{(1)}=[0,2p]$.
\item[ for $s = 2$: ]
$H_k^{(2)} = \displaystyle \frac{H^{(1)}_{k-1} + H^{(1)}_{k}}{2}$
\- and $\displaystyle \Phi^{(2)}$ admits an approximation of degree
$p+2$ and $\supp\displaystyle \Phi^{(1)}=[0,2p+1]$. In addition, according to
the step $s=1$, we can rewrite: $$\displaystyle H_k^{(2)} =
\displaystyle \frac{H^{(0)}_k+ H^{(0)}_{k-1} + H^{(0)}_{k-1} +
H^{(0)}_{k-2}}{2^2}.$$
\item[] \noindent The result for $s>2$ is obtained by replacing $\Phi^{(s)},\,\tilde\Phi^{(s)}$ by $\Phi^{(s+1)},\,\tilde \Phi^{(s+1)}$
and applying succesively the result for $s=1$.
\end{description}
The dual discrete filters are obtained by substituting $h_k$ by $\tilde{H}_k$ and
$\tilde{h}_k$ by $H_k$ with Property (\ref{regular}).$\square$\\
We conclude with some representation of wavelets in the biorthogonal and orthonormal cases:
\begin{figure}[H]
\centering{
\includegraphics[scale=0.25]{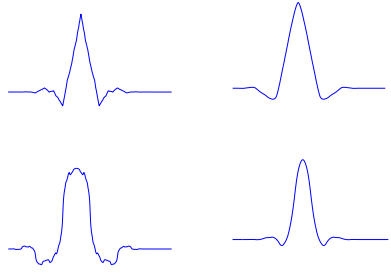}
}
\caption{Top left $\phi$, top  right $\tilde{\phi}$,
bottom left $\Phi^{(1)}$, bottom right $\tilde{\Phi}^{(1)}$ where
$\phi,\tilde\phi$ are
    compactly supported B-spline 4 4}\label{Fig1}
\end{figure}
\begin{figure}[H]
\centering{
\includegraphics[scale=0.25]{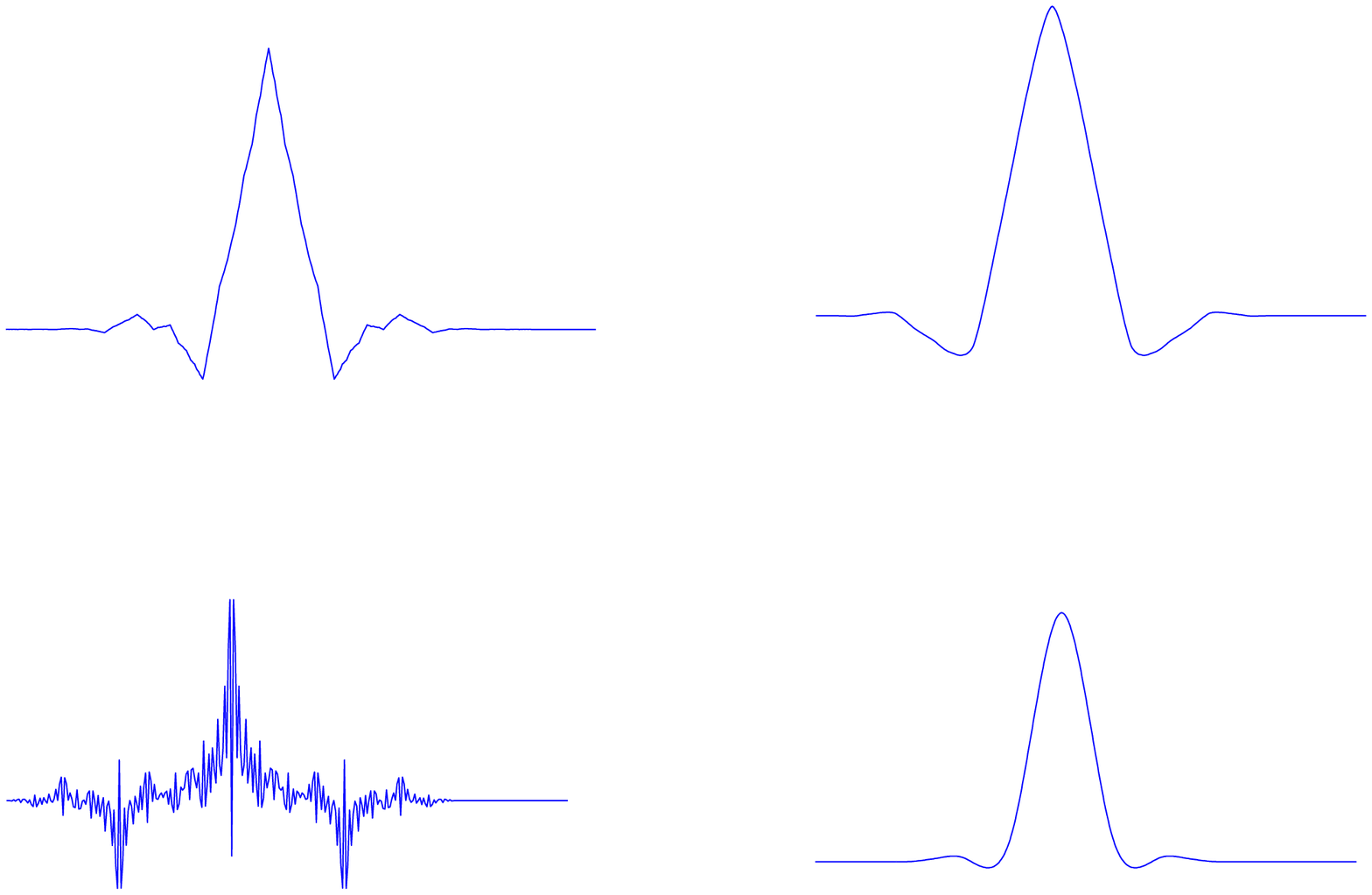}
}
\caption{Top left $\phi$, top  right $\tilde{\phi}$,
bottom left $\Phi^{(2)}$, bottom right $\tilde{\Phi}^{(2)}$ where
$\phi,\tilde\phi$ are
    compactly supported B-spline 4 4 }\label{Fig2}
\end{figure}
\begin{figure}[H]
\centering{
\includegraphics[scale=0.25]{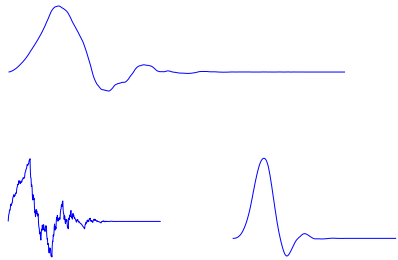}
}
\caption{Top left $\phi$, bottom left $\Phi^{(1)}$,
bottom right $\tilde{\Phi}^{(1)}$ where $\phi$ is the Daubechies
Wavelet with $4$ vanishing moments}\label{Fig4}
\end{figure}
\section{Conclusion}
The constructed wavelet $\tilde\Phi^{(s)}$ is more regular  than $\Phi^{(s)}$ but has
less vanishing moments. We have gained vanishing moments for one and  regularity for the other but
in counterpart we have increased the support of the scaling function but still compact. This operation preserves
the parity of the initial wavelet and if initial wavelets are curl free then  so are the
constructed wavelets. Moreover,  this construction is a way to recover  all boxed splines.
We have also a similar way to construct biorthogonal  wavelets from an orthonormal initial family defining $\widetilde\Phi^{(s)} =
\left(T^{\star}\right)^{(s)} \phi$.
\bibliographystyle{plain}

\end{document}